\documentclass[10pt]{article}
\usepackage{color,graphicx}
\newtheorem{theorem}{Theorem} 
\newtheorem{lema}[theorem]{Lemma}
\newtheorem{obs}[theorem]{Observation}

\newtheorem{ques}[theorem]{Question}

\newcommand{\Z}{\hbox{\bf Z}}
\newcommand{\F}{\hbox{\bf F}}

\newcommand{\beeq}{\begin{eqnarray*}}
\newcommand{\eneq}{\end{eqnarray*}}
\newcommand{\proof}{\noindent {\it Proof.\hspace{2mm}}}

\newcommand{\qfd}{\hfill $\fbox{}$\vspace{2mm}}
\def\newpic#1{%
\def\emline##1##2##3##4##5##6{%
\put(##1,##2){\special{em:point #1##3}}%
\put(##4,##5){\special{em:point #1##6}}%
\special{em:line #1##3,#1##6}}}
\newpic{}
\def\emline#1#2#3#4#5#6{%
\put(#1,#2){\special{em:moveto}}%
\put(#4,#5){\special{em:lineto}}}
\def\newpic#1{}
\newcommand\ZZ{{{\rm Z}\kern-.28em{\rm Z}}}

\title{Worst-case efficient dominating sets in digraphs}
\author{Italo J. Dejter
\\ University of Puerto Rico \\ Rio Piedras, PR 00936-8377 \\ italo.dejter@gmail.com
}
\date{}

\begin{document}
\maketitle

\begin{abstract}
\noindent Let $1\le n\in\Z$. {\it Worst-case efficient dominating
sets in digraphs} are conceived so that their presence in certain
strong digraphs $\vec{ST}_n$ corresponds to that of efficient
dominating sets in star graphs $ST_n$:
The fact that the star graphs $ST_n$ form a so-called dense
segmental neighborly E-chain is reflected in a
corresponding fact for the digraphs $\vec{ST}_n$. Related
chains of graphs and open problems are presented as well.
\end{abstract}

\noindent{\bf Keywords:} Cayley graph; star graph; digraph;
efficient dominating set

\section{Introduction}

\noindent In this work, it is shown that worst-case efficient
dominating sets $S$ in digraphs, introduced in the next paragraphs,
play a role in oriented Cayley graph variants $\vec{ST}_n$ of the
star graphs $ST_n$ \cite{akers} that adapts the role played by
efficient dominating sets \cite{arum1,Hay} in the $ST_n$ ($1\le
n\in\Z$). That the $ST_n$ form a so-called dense segmental
neighborly E-chain \cite{io} is then reflected (Section 5) in a
corres\-ponding property for the $\vec{ST}_n$\,, which are formally
defined in Section 3 and henceforth referred to as the {\it star
digraphs}. A non-dense non-neighborly version of this reflection
goes from perfect codes in binary Hamming cubes to worst-case
perfect codes (another name for worst-case efficient dominating
sets) in oriented ternary Hamming cubes. In Section 6, some comments
and open problems on hamiltonicity and traceability of these star
digraphs and on related concepts of {\it pancake} and {\it
binary-star digraphs} are also presented.\bigskip

\noindent Let $D$ be a digraph. A vertex $v$ in $D$ is said to be a
{\it source}, (resp. a {\it sink}), in $D$ if its indegree
$\partial^-(v)$ is null, (resp. positive), and its outdegree
$\partial^+(v)$ is positive, (resp. null). A vertex subset $S$ of
$D$ is said to be {\it worst-case stable}, or {\it $\pm$stable}, if
$\min\{\partial^-(v),\partial^+(v)\}=0$, for every vertex $v$ in the
induced subdigraph $D[S]$, or equivalently: if each vertex in $D[S]$
is either a source or a sink or an isolated vertex. In this case,
$D[S]$ is a directed graph with no directed cycles, called a {\it
directed acyclic graph}, or {\it dag}. On the other hand, a vertex
subset $S$ in $D$ is said to be {\it stable} if it is stable in the
underlying undirected graph of $D$. Clearly, every stable $S$ in $D$ is $\pm$stable, but the converse is not true in general.\bigskip

\noindent Given a vertex $v$ in $D$, if there exists an arc $(u,v)$,
(resp. $(v,u)$), in $D$, then we say that $v$ is $(+)$dominated,
(resp. $(-)$dominated), by $u$, and that $v$ is a $(+)$neighbor,
(resp. $(-)$neighbor), of $u$. Given a subset $S$ of vertices of
$D$, if each vertex $v$ in $D-S$ is $(+)$dominated by a vertex $u$
 in $S$ and $(-)$dominated by a vertex $w$ in $S$, then we say that
$S$ is a {\it $\pm$dominating set} in $D$. If $u$ and $w$ are
unique, for each vertex $v$ in $D$, then $S$ is said to be {\it
perfect}. If $D$ is an oriented simple graph, then in the previous
sentence it is clear that $u\ne w$. A vertex subset in $D$ is said
to be a {\it worst-case efficient dominating set} if it is both
perfect $\pm$dominating and $\pm$stable. The {\it worst-case
domination number} $\gamma_\pm(D)$ of $D$ is the minimum cardinality
of a worst-case efficient dominating set in $D$.\bigskip

\noindent Given a vertex set $S$ in $D$, let $N^-(S)$, (resp.
$N^+(S)$), be the subset of vertices $u$ in $V(D)\setminus S$ such
that $(u,v)$, (resp. $(v,u)$), is an arc in $D$, for some vertex $v$
 in $S$. A vertex subset $S$ in $D$ is said to be {\it cuneiform} if
$N^+(S)\cup N^-(S)$ is the disjoint union of two stable vertex
subsets in $D$, as indicated, and there is a bijective
correspondence $\rho:S\rightarrow K$ such that $v$ and $\rho(v)$
induce a directed triangle $\vec{\Delta}_v$\,, for each $v\in S$,
where $K$ is a disjoint union of $|S|$ digraphs $\vec{P}_2$ in $D$
consisting each of a single arc from $N^+(S)$ to $N^-(S)$. A
worst-case efficient dominating set $S$ in $D$ is said to be an {\it
E$_\pm$-set} if it cuneiform. A subdigraph in $D$ induced by an
E$_\pm$-set is said to be an {\it E$_\pm$-subdigraph}.\bigskip

\noindent In undirected graphs, E-sets \cite{io} correspond to
perfect ($1$-error-correcting) codes \cite{bi,Kr}. A version for
E$_\pm$-sets of the sphere-packing condition for E-sets in \cite{io}
is given as follows: If an oriented simple graph $D$ has the same
number $r$ of $(+)$neighbors as it has of $(-)$neighbors at every
vertex so that the outdegree and the indegree of every vertex are
both equal to $r$, then
\begin{eqnarray}|V(D)|=(2+r)|S|/2,\end{eqnarray}
for every E$_\pm$-set $S$ in $D$, where the factor $(2+r)$ accounts
for each of the $|S|/2$ sources and each of the $|S|/2$ sinks of
$S$, and for the $r$ $(+)$neighbors of each such source, or
alternatively the $r$ $(-)$neighbors of each such sink. Clearly, the
{\it oriented sphere-packing condition} (1) is a necessary condition
for $S$ to be an E$_\pm$-set in $D$. (Compare this, and other
concepts defined in Section 2, with \cite{io}).

\section{E$_\pm$-chains}

\noindent Inspired by the concept of E-chain for undirected graphs
in \cite{io}, a countable family of oriented simple graphs disposed
as an increasing chain by containment

\begin{eqnarray}{\mathcal D}=\{ D_1\subset D_2\subset\ldots\subset
D_n\subset D_{n+1}\subset\ldots\}\end{eqnarray} is said to be an
{\it E$_\pm$-chain} if every $D_n$ is an induced subdigraph of
$D_{n+1}$ and each $D_n$ contains an E$_\pm$-set $S_n$. For oriented
simple graphs $D$ and $D'$\,, a one-to-one digraph map
$\zeta:D\rightarrow D'$ is an {\it inclusive map} if $\zeta(D)$ is
an induced subdigraph in $D'$. This is clearly an {\it
orientation-preserving map}, or {\it $(+)$map}, but we also consider
{\it orientation-reversing maps}, or {\it $(-)$maps},
$\zeta:D\rightarrow D'$ and say that one such map is {\em inclusive}
if $\zeta(D)$ is an induced subdigraph in $D'$\,, even though
corresponding arcs in $D$ and in $\zeta(D)$ are oppositely oriented
in this case.\bigskip

\noindent Let
\begin{eqnarray}\kappa_n:D_n\rightarrow D_{n+1}\end{eqnarray}
stand for the inclusive map of $D_n$ into $D_{n+1}$ induced by
${\mathcal D}$, where $n\ge 1$. If $V(\kappa_n(D_n))$ is cuneiform,
for every $n\ge 1$, then we say that the E$_\pm$-chain ${\mathcal
D}$ is a {\it neighborly} E$_\pm$-chain.\bigskip

\noindent If there exists an inclusive map
\begin{eqnarray}\zeta_n: D_n \rightarrow D_{n+1}\end{eqnarray} such
that $\zeta_n(S_n)\subset S_{n+1}$\,, for each $n\ge 1$, then we say
that the E$_\pm$-chain ${\mathcal D}$ is {\it inclusive}, (where
each $\zeta_n$ is either a $(+)$map or a $(-)$map). Notice that an
inclusive neighborly E$_\pm$-chain has $\kappa_n\neq\zeta_n$\,, for
every integer $n\ge 1$. A particular case of inclusive E$_\pm$-chain
$\mathcal D$ is one in which $S_{n+1}$ has a partition into images
$\zeta_n^{(k)}(S_n)$ of $S_n$ through respective inclusive maps
$\zeta_n^{(k)}$\,, where $k$ varies on a suitable finite indexing
set. In such a case, the E$_\pm$-chain $\mathcal D$ is said to be
{\it segmental}.\bigskip

\noindent An E$_\pm$-chain ${\mathcal D}$ of oriented simple graphs
that have the same number $r$ of $(+)$neighbors as it has of
$(-)$neighbors at every vertex so that the outdegree and the
indegree of every vertex are both equal to $r$ is said to be {\it
dense} if $|S_n|/|V(D_n)|=2/(n+1)$\,, for each $n\ge 1$, in
accordance with the modified sphere-packing condition (1). It can be
seen \cite{io} that the star graphs $ST_n$ form a dense segmental
neighborly E-chain, while the Hamming cubes $F(2,2^n-1)$ form a
segmental E-chain which is neither neighborly nor dense, where $n\ge
1$. An example of a dense segmental neighborly E$_\pm$-chain is
given by the star digraphs $D_n=\vec{ST}_{n+1}$\,, to be treated in
Sections 3 to 5 below. On the other hand, we note that the ternary
Hamming cubes $F(3,\frac{3^n-1}{2})$\,, with edge orientations
induced by the order $0<1<2$ in the 3-element field $\F_3=\Z/3\Z$,
form a segmental E$_\pm$-chain which is neither neighborly nor
dense. These cubes considered undirected constitute a segmental
E-chain which is neither neighborly nor dense. However, in their
oriented version, their E$_\pm$-sets may be called now {\it
worst-case perfect codes} (another name for E$_\pm$-sets). The
E$_\pm$-sets in the star digraphs $\vec{ST}_n$ taken undirected are
not E-sets, though.

\section{Star digraphs}

\noindent Let $n\ge 1$. The {\it star graph} $ST_n$ is the Cayley
graph of the group $Sym_n$ of symmetries on the set
$\{0,1,\ldots,n-1\}$ with respect to the generating set formed by
the transpositions $(0\;i)$, where $i\in\{1,\ldots,n-1\}$. While
$ST_2=K_2$ and $ST_3$ is a 6-cycle graph, the undirected graph
induced by the curved arcs in Figure 1 below (with symbol 4 omitted
in each vertex) shows a cutout of $ST_4$ embedded into a toroid $T$
obtained by identification of the 3 pairs of opposite sides in the
external intermittent hexagon. In the figures of this section,
elements of $Sym_n$\,, or of its alternating subgroup $Alt_n$\,, are
represented by $n$-tuples $x_0x_1\ldots x_{n-1}$ corresponding to
respective permutations ${0\;\;1\;\ldots\ldots\; n-1 \choose
x_0x_1\ldots\ldots \;x_{n-1}}$ in $Sym_n$ or $Alt_n$.\bigskip

\begin{figure}[htp]
\vspace*{-3mm} \hspace*{7.0mm}
\includegraphics[scale=0.35]{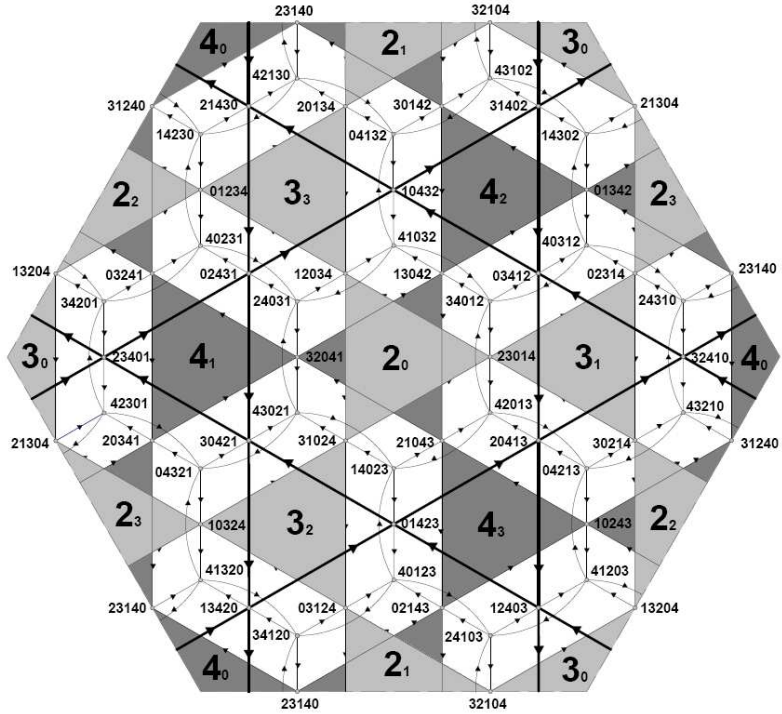}\\
\vspace*{-6mm} \caption{Toroidal cutout of $\vec{ST}_5$ stressing
subdigraphs $\vec{ST}_4^{i,j}$}
\end{figure}

\noindent The {\it star digraph} $\vec{ST}_n$ is the Cayley graph of
$Alt_n$ with respect to a generating set formed by the permutations
$(0\;1\;i)=(0\;1)(1\; i)$, where $i\in\{2,\ldots,n-1\}\}$. This is
an oriented simple graph, so that it does not have cycles of length
2. While $\vec{ST}_2=K_1$ and $\vec{ST}_3$ is a directed triangle,
$\vec{ST}_4$ is an edge-oriented cuboctahedron, as depicted in any
of the 4 instances in Figure 2. This and Figures 1 and 3, showing
different features in $\vec{ST}_5$\,, are presented
subsequently.\bigskip

\noindent The digraph $\vec{ST}_5$\,, with its 60 vertices
corresponding to (and identified with) the 60 even permutation on
the set $\{0,1,2,3,4\}$, is the edge-disjoint union of 4 induced
subdigraphs embedded into the toroid $T$ cut out in Figure 1. Shown
in this figure, these 4 toroidal subdigraphs are: {\bf(i)} the
subdigraph spanned by the 36 oriented 3-cycles $(v,w,u)$ with two
short straight contiguously colinear arcs $(v,w),(w,u)$ and a longer
returning curved arc $(u,v)$ having tail $u$, (resp. head $v$), with
first, (resp. second), entry equal to 4; {\bf(ii)} the subdigraph
$\vec{ST}_4^{4,4}$ on 12 vertices spanned by the light-gray
triangles, whose vertices have their fifth entry equal to 4;
{\bf(iii)} the subdigraph $\vec{ST}_4^{4,3}$ on 12 vertices spanned
by the dark-gray triangles (some partially hidden by the light-gray
triangles), whose vertices have their fourth entry equal to 4;
{\bf(iv)} the subdigraph $\vec{ST}_4^{4,2}$ on 12 vertices spanned
by the triangles with bold-traced arcs, whose vertices have their
third entry equal to 4.\bigskip

\begin{figure}[htp]
\vspace*{-3mm} \hspace*{0.5mm}
\includegraphics[scale=0.4]{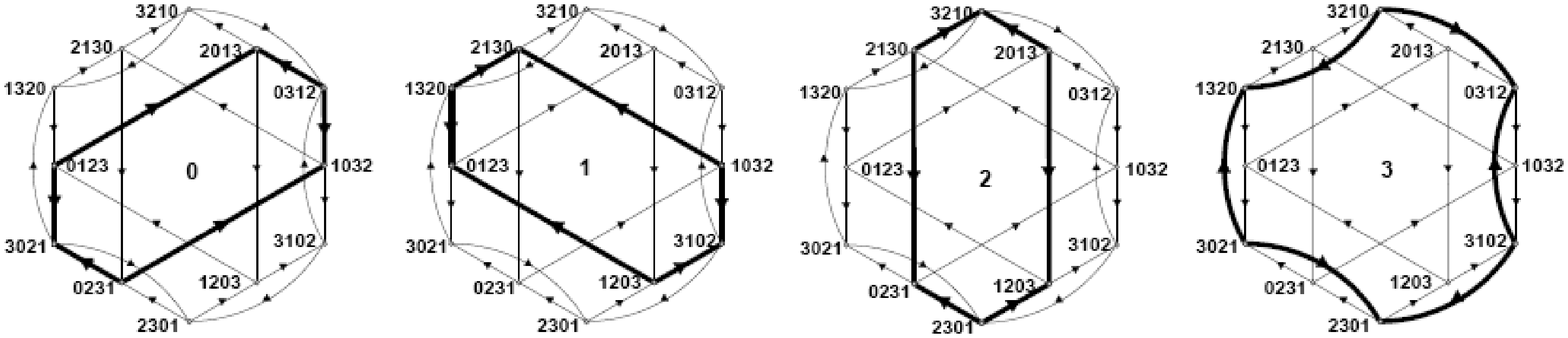}\\
\vspace*{-6mm} \caption{Representations of $\vec{ST}_4$ stressing
$ST_3^0$\,, $ST_3^1$\,, $ST_3^2$ and $ST_3^3$}
\end{figure}

\noindent Let us indicate the first to fifth entries of the vertices
of $\vec{ST}_5$ as 0th to 4th entries, respectively. The 24-vertex
subdigraph $ST_4^4$ induced by the 36 curved arcs in item (i) above
is an induced copy of $ST_4$ with its edges oriented from vertices
with 0th entry equal to 4 to vertices with 1st entry equal to 4. On
the other hand, the subdigraphs $\vec{ST}_4^{4,j}$ above, (for
$j=2,3,4$\,, corresponding to items (iv), (iii), (ii),
respectively), form 3 induced copies of $\vec{ST}_4$ in
$\vec{ST}_5$. Moreover, there are 15 induced copies of $\vec{ST}_4$
in $\vec{ST}_5$: Apart from the digraphs $\vec{ST}_4^{4,j}$ above,
the symbols $j_i=2_0,\ldots, 4_3$ in Figure 1 denote the other 12
induced copies of $\vec{ST}_4$ in $\vec{ST}_5$\,, each with 6 curved
arcs forming a 6-cycle dag, (with each two contiguous arcs having
opposite orientations), and the symbol $j_i$ displayed at its
center. In order to maintain the notation of the 3 initially
presented copies of $\vec{ST}_4$ in $\vec{ST}_5$\,, we denote these
12 new copies in their order of presentation above as follows:
$$\vec{ST}_4^{i,j}=\vec{ST}_4^{0,2}, \ldots, \vec{ST}_4^{3,4}.$$
Each of these subdigraphs is induced by all the vertices with $j$th
entry equal to $i$\,, where $i\in\{0,1,2,3\}$ and $j\in\{2,3,4\}$.
These subdigraphs are the images of corresponding maps
$\zeta_4^{i,j}=\xi_4^{i,j}$ whose definition is completed in Section
4.\bigskip

\noindent Figure 2 shows 4 representations of $\vec{ST}_4$. In each,
the edges of a distinctive copy of the star graph $ST_3$ are shown
in bold trace, with each edge oriented from a vertex with 0th entry
equal to $i$ to a vertex with 1st entry equal to $i$\,, where
$i=0,1,2,3$\,, respectively. This symbol $i$ is shown at the center
of the corresponding representation. Accordingly, we denote these
copies by $ST_3^0$\,, $ST_3^1$\,, $ST_3^2$ and $ST_3^3$\,,
respectively. Observe that each such oriented $ST_3^i$ is an
E$_\pm$-set so that $\vec{ST}_4-ST_3^i$ is the disjoint union of two
copies of $\vec{ST}_3$ (directed triangles) that we denote
$\vec{ST}_3^{i,3}$ and $\vec{ST}_3^{i,2}$ according to whether the
3rd or the 2nd entry of its vertices is equal to
$i=0,1,2,3$.\bigskip

\begin{figure}[htp]
\vspace*{-3mm} \hspace*{7.0mm}
\includegraphics[scale=0.35]{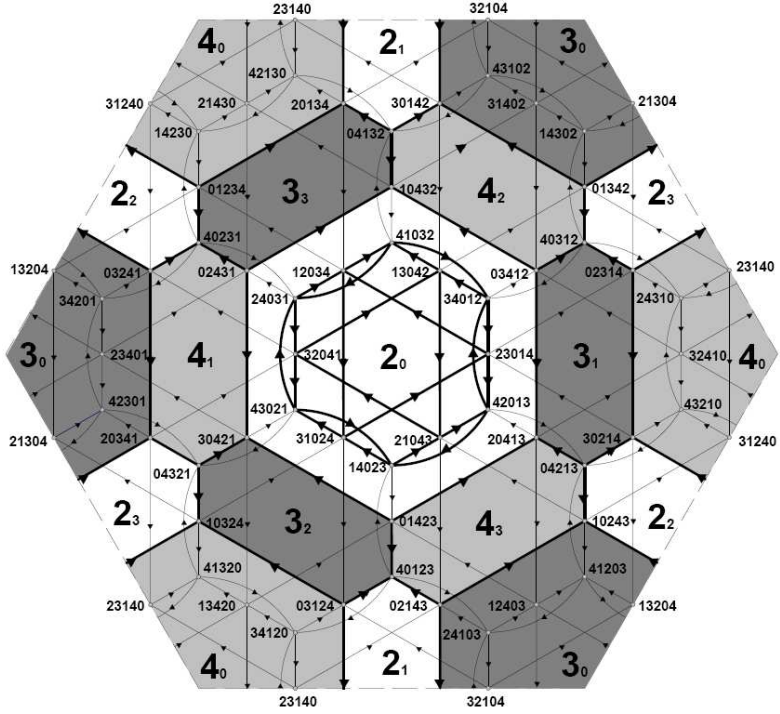}\\
\vspace*{-6mm} \caption{Toroidal cutout of $\vec{ST}_5$ stressing
$\vec{ST}_4^{0,2}=2_0$ and its E$_\pm$-set $ST_4^0$}
\end{figure}

\noindent Recall that the subdigraph $ST_4^4$ of $\vec{ST}_5$ induced
by the set of curved arcs in Figure 1 constitutes an induced copy of
an (undirected) $ST_4$ in $\vec{ST}_5$ with its edges considered as
arcs of $\vec{ST}_5$ by orienting them as shown, that is from those
vertices with their 3 departing curved arcs being counterclockwise
concave to those vertices with their 3 arriving curved arcs being
clockwise concave. They constitute respectively vertex parts
$V_0(ST_4^4)$ and $V_1(ST_4^4)$ of $ST_4^4$ considered as a
bipartite graph, where $V_0(ST_4^4)$\,, (resp. $V_1(ST_4^4)$), is
formed by those vertices of $\vec{ST}_5$ with 0th, (resp. 1st),
entry equal to 4.\bigskip

\noindent In terms of our desired chain, as in (2), that we want
made up of star digraphs

\begin{eqnarray}{\mathcal D}=\{\vec{ST}_2\subset\vec{ST}_3\subset
\ldots\subset\vec{ST}_n\subset
\vec{ST}_{n+1}\subset\ldots\},\end{eqnarray} where
$D_1=\vec{ST}_2$\,, $D_2=\vec{ST}_3$\,, $\ldots$\,,
$D_n=\vec{ST}_{n+1}$\,, etc., the inclusions $\vec{ST}_n\subset
\vec{ST}_{n+1}$ must be given via corresponding maps
$\kappa_{n-1}:\vec{ST}_n\rightarrow \vec{ST}_{n+1}$. This fits for
$n=4$ as $\vec{ST}_4^{4,4}=\kappa_3(\vec{ST}_4)$\,, represented by
the graph spanned via the light-gray triangles in Figure 1.\bigskip

\noindent Figure 3\, is as Figure 1, but showing centrally the
induced copy $2_0=\vec{ST}_4^{0,2}$ of $\vec{ST}_4$ in $\vec{ST}_5$
with its arcs in bold trace to stress that the neighbors of its
vertices induce an oriented\, $\pm$stable copy of $ST_4$ whose 12
6-cycle dags (with contiguous arcs oppositely oriented) are shown
delimiting 4 light-gray, 4 dark-gray and 4 white colored regions.
Thus, such a copy of $ST_4$ contains 12 oriented $\pm$stable copies
of $ST_3$\,, namely those 12 6-cycle dags. Each of these copies may
be used to define a map $\zeta$ as in (4) to prove that the chain
$\mathcal D$ in (5) above is segmental.\bigskip

\noindent Observe from Figures 1, 2 and 3 that $\vec{ST}_4$\,,
(resp. $\vec{ST}_5$), has 2 and 6, (resp. 3 and 12), maps in the
nature of $\kappa_i$ and $\zeta_i^{(k)}$\,, respectively, as defined
in Section 2, ($i=2,3$). We illustrate this after formalizing some
definitions.

\section{Formalizing definitions}

\noindent For $1<n\in\Z$\,, the Greek letters $\kappa$ and $\zeta$
used since (3) and (4) above will be consolidated as letter $\zeta$
in defining formally some useful graph maps, below:
$$\zeta_n^{i,n}=\kappa_n^{i,n}:\vec{ST}_n\rightarrow\vec{ST}_{n+1}
\;\;\;\mbox{ and }\;\;\;
\zeta_n^{i,j}:\vec{ST}_n\rightarrow\vec{ST}_{n+1},$$ where $1<j<
n$\,,\, $0\le i\le n$ and the $\kappa$-notation as in (3) (resp.
$\zeta$-notation as in (4)) is associated with the neighborly (resp.
segmental) E$_\pm$-chains defined in Section 2, so that graph
inclusions $\kappa_{n-1}$ as in (3) will now be denoted
$\zeta_n^{n,n}$ (or $\kappa_n^{n,n}$) and so that we could also
denote
$$\vec{ST}_n^{i,j}=\zeta_n^{i,j}(\vec{ST}_n)\subset\vec{ST}_{n+1},$$
accompanying the notation of the examples of Section 3 above.\bigskip

\noindent The maps $\zeta_n^{i,j}$ are defined as follows: For each
$i\in I_{n+1}=\{0,1,\ldots,n\}$\,, let
$$\phi^i:I_n\rightarrow I_{n+1}$$ be given by
$$\phi^i(x) =
\left\{
    \begin{array}{ll}
        x  & \mbox{, if } 0\le x<i;  \\
        x+1  & \mbox{, if } i\le x<n.
    \end{array}
\right.$$

\noindent We denote $\phi^i(x)=x^i$. For $0\le i\le n$\, and\, $2\le
j\le n$\,, we define
$$\zeta_n^{i,j}(x_0x_1\ldots x_{n-1})=\psi^{i,j}(x_0^i,x_1^i)\ldots x_{j-1}^ii
x_j^i\ldots x_{n-1}^i,$$ where
$$\psi^{i,j}(xy) =
\left\{
    \begin{array}{ll}
        xy  & \mbox{, if } (n-i+j) \equiv 0 \mbox{ (mod 2)}; \\
        yx  & \mbox{, if } (n-i+j) \equiv 1 \mbox{ (mod 2)}.
    \end{array}
\right.$$

\noindent  To continue the examples from Section 3, observe that:
{\bf(a)} the $(+)$map $\kappa_2=\kappa_3^{3,3}=\zeta_3^{3,3}$ sends
$x_0x_1x_2\in V(\vec{ST}_3)$ onto $x_0x_1x_23$\,, for each
$\{x_0,x_1,x_2\}=\{0,1,2\}$; {\bf(b)} the $(+)$map
$\kappa_3=\zeta_4^{4,4}$ sends the vertices $x_0x_1x_2x_3$ of
$\vec{ST}_3$ onto the corresponding vertices $x_0x_1x_2x_34$ of
$\vec{ST}_4$\,, for each $\{x_0,x_1,x_2,x_3\}=\{0,1,2,3\}$; {\bf(c)}
the $(-)$map $\kappa_3^{3,2}=\zeta_3^{3,2}$ sends the vertices
$x_0x_1x_2$ of $\vec{ST}_3$ onto the corresponding vertices
$x_1x_03x_2$ of $\vec{ST}_4$; {\bf(d)} the $(-)$map
$\kappa_4^{4,3}=\zeta_4^{4,3}$ sends the vertices $x_0x_1x_2x_3$ of
$\vec{ST}_4$ onto the corresponding vertices $x_1x_0x_24x_3$ of
$\vec{ST}_5$ and {\bf(e)} the $(+)$map
$\kappa_4^{4,2}=\zeta_4^{4,2}$ sends the vertices $x_0x_1x_2x_3$ of
$\vec{ST}_4$ onto the corresponding vertices $x_0x_14x_2x_3$ of
$\vec{ST}_5$. Note that the image of $\kappa_n^{n,i}$ is
$\vec{ST}_n^{n,i}$ for every $n>1$\,, starting at $\kappa_1(01)=012$
for
$\kappa_1=\kappa_2^{2,2}:\vec{ST}_2\rightarrow\vec{ST}_3$.\bigskip

\noindent The images of maps $\zeta_1^{(k)}$ in Figure 2 are those
directed triangles having one curved arc: There are 6 such images.
The images of maps $\zeta_2^{(k)}$ in the digraph $\vec{ST}_5$ of
Figures 1 and 3 look like any of the 4 representations of
$\vec{ST}_4$ in Figure 2 (disregarding the thickness of their
edges): There are 12 such images. In general, for $2\le i<n$ and
$0\le j<n$ it holds that $\vec{ST}_n^{i,j}$ is the image of such a
map $\zeta_n^{i,j}$ (or $\zeta_n^{i,j}$), which is either a $(+)$map
or a $(-)$map.\bigskip

\noindent To be applied in the proof of Theorem 6 below, the
E$_\pm$-set $ST_4^0$ in Figure 3 admits 3 different partitions into
4 6-cycle dags, each composed by alternating sources and sinks: One
partition with 6-cycle-dag interiors in light-gray color, another
one in dark-gray color and the third one with white interiors. A
table of the involved quadruples of 6-cycle dags of $\vec{ST}_5$
follows, with the 6-cycle dag
$ST_3^0=(0123>2013<0312<1032<0231>3021<)$ (having $>$ and $<$
standing for forward and backward arcs between contiguous vertices,
the sixth vertex taken contiguous to the first vertex) sent via the
graph maps $\zeta_4^{i,j}$ onto the corresponding vertices of
$\vec{ST}_5$ in the quadruples, and with $\ge$, (resp. $\le$),
standing for coincidence in case $\zeta_4^{i,j}$ is a $(+)$map,
(resp. $(-)$map):
$$\begin{array}{l}
^{\{\zeta_4^{1,4}(
ST_3^0)\,\le\,(20341<03241>40231<02431>30421<04321>)\subset\vec{ST}_4^{1,4},}_{\;\;\zeta_4^{2,4}(
ST_3^0)\,\ge\,(01342>30142<04132>10432<03412>40312<)\subset\vec{ST}_4^{2,4},}\vspace*{1mm}\\^{\;\;\zeta_4^{3,4}(
ST_3^0)\,\le\,(10243<02143>40123<01423>20413<04213>)\subset\vec{ST}_4^{3,4},}_{\;\;\zeta_4^{4,4}(
ST_3^0)\,\ge\,(01234>20134<03124>10324<02314>30214<)\subset\vec{ST}_4^{4,4}\},}\vspace*{2mm}\\^{\{\zeta_4^{1,3}(
ST_3^0)\,\ge\,(02314>30214<04213>20413<03412>40312<)\subset\vec{ST}_4^{1,3},}_{\;\;\zeta_4^{2,3}(
ST_3^0)\,\le\,(10324<03124>40123<01423>30421<04321>)\subset\vec{ST}_4^{2,3},}\vspace*{1mm}\\^{\;\;\zeta_4^{3,3}(
ST_3^0)\,\ge\,(01234>20134<04132>10432<02431>40231<)\subset\vec{ST}_4^{3,3},}_{\;\;\zeta_4^{4,3}(
ST_3^0)\,\le\,(10243<02143>30142<01342>20341<03241>)\subset\vec{ST}_4^{4,3}\},}\vspace*{2mm}\\^{\{\zeta_4^{1,2}(
ST_3^0)\,\le\,(20134<03124>40123<02143>30142<04132>)\subset\vec{ST}_4^{1,2},}_{\;\;\zeta_4^{2,2}(
ST_3^0)\,\ge\,(01234>30214<04213>10243<03241>40231<)\subset\vec{ST}_4^{2,2},}\vspace*{1mm}\\^{\;\;\zeta_4^{3,2}(
ST_3^0)\,\le\,(10324<02314>40312<01342>20341<04321>)\subset\vec{ST}_4^{3,2},}_{\;\;\zeta_4^{4,2}(
ST_3^0)\,\ge\,(01423>20413<03412>10432<02431>30421<)\subset\vec{ST}_4^{4,2}\},}\\
\end{array}$$
namely $\{\zeta_4^{i,j}(ST_3^0)\,|\,i=1,2,3,4\}$, for $j=4,3,2$,
resp., where the composing vertices of each target 6-cycle dag
correspond orderly with those of the source 6-cycle dag. Each target
6-cycle dag here is, for $n=3$, the induced digraph of a E$_\pm$-set
$S_n$ in a corresponding subdigraph $D_n=\vec{ST}_{n+1}^{i,j}$,
required for every $n\ge 1$ to insure that $\mathcal D$ is
inclusive, as defined in the paragraph containing display (4). If
partitions as in that paragraph are obtained for every $n\ge 1$\,,
as they were for $n=3$ above, then $\mathcal D$ is segmental.

\section{Main results}

\begin{obs} Let $1<n\in\Z$ and let\, $0\le
i\le n$\, and\, $2\le j\le n$. If $(n-i+j)$ is even, $($resp.
odd$)$\,, then the map $\zeta_n^{i,j}$ is a $(+)$map, $($resp.
$(-)$map$)$.
\end{obs}

\proof The auxiliary assignment $\psi^{i,j}$ above is devised so as
to allow keeping the parity of the target permutations of the map
$\zeta_n^{i,j}$. A side effect of this is that, depending on the
parity of the quantity $(n-i+j)$\,, the maps $\zeta_n^{i,j}$ are
compelled to be alternatively $(+)$maps and $(-)$maps. To illustrate
this point, observe that for $n=3$ the maps $\zeta_n^{i,j}$ are:
\begin{enumerate}
\item the $(+)$map $\zeta_2^{0,2}:\vec{ST}_2\rightarrow\vec{ST}_3$ given by
$\zeta_2^{0,2}(01)=(120)$;
\item the $(-)$map $\zeta_2^{1,2}:\vec{ST}_2\rightarrow\vec{ST}_3$ given by
$\zeta_2^{1,2}(01)=(201)$;
\item the $(+)$map $\zeta_2^{2,2}:\vec{ST}_2\rightarrow\vec{ST}_3$ given by
$\zeta_2^{2,2}(01)=(012)$.
\end{enumerate}
The neighbors of each of these copies $\vec{ST}_2^{i,2}$ of
$\vec{ST_2}$ induce a corresponding copy of $ST_2$ with its (only)
edge oriented from the vertex with 0th entry equal to $i$ into the
vertex with 1st entry equal to $i$\,, for $i=0,1,2$\,, respectively.
\qfd

\begin{lema}
For $n>1$\,, there are $n^2-1$ copies of the star subdigraph
$\vec{ST}_n$ in the star digraph $\vec{ST}_{n+1}$. These copies are
the images of $\vec{ST}_n$ under the $\lceil\frac{n^2-1}{2}\rceil$\,
$(+)$maps and $\lfloor\frac{n^2-1}{2}\rfloor$\, $(-)$maps
$\zeta_n^{i,j}$\,, where $0\le i\le n$ and $2\le j\le n$.
\end{lema}

\proof Notice that $n^2-1=(n+1)(n-1)$ is the number of digraph maps
$\zeta_n^{i,j}$ defined before Observation 1. Of them, $\lceil
\frac{n^2-1}{2}\rceil$ are $(+)$maps and
$\lfloor\frac{n^2-1}{2}\rfloor$ are $(-)$maps, depending by
Observation 1 on the condition that $(n-i+j)$ be even or odd,
respectively. Also, the image of each $\zeta_n^{i,j}$ is the
corresponding copy $\vec{ST}_n^{i,j}$ of $\vec{ST}_n$ in
$\vec{ST}_{n+1}$. This yields the statement. \qfd

\begin{theorem} The worst-case domination number of
$\vec{ST}_{n+1}$ is $\gamma(\vec{ST}_{n+1})=n!$
\end{theorem}

\proof Since $|V(\vec{ST}_{n+1})|=\frac{(n+1)!}{2}$ and
$|V(\vec{ST}_n^{i,j})|=\frac{n!}{2}$\,, for each
$i\in\{0,\ldots,n\}$ and $j\in\{2,\ldots,n\}$\,, we have that
$N^+(V(\vec{ST}_n^{i,j}))\cup N^-(V(\vec{ST}_n^{i,j}))$ is an
E$_\pm$-set in $\vec{ST}_{n+1}$ inducing the E$_\pm$-digraph
$ST_n^j$\,, for which $|V(ST_n^j)|=n!$\, Indeed, as exemplified in
Figure 3, $S=N^+(V(\vec{ST}_n^{i,j}))\cup N^-(V(\vec{ST}_n^{i,j}))$
is {\bf(a)} $\pm$stable, for its induced subdigraph
$\vec{ST}_{n+1}[S]$ is composed by sources and sinks; {\bf(b)}
$\pm$dominating, for each vertex $v$ in $\vec{ST}_{n+1}-S$ is
$(+)$dominated by a vertex $u$ in $S$ and $(-)$dominating by a
vertex $w$ in $S$; and {\bf(c)} perfect, for the vertices $u$ and
$w$ in (b) above are unique, for each vertex $v$ of
$\vec{ST}_{n+1}-S$; {\bf(d)} an E$_\pm$-set, for the definition of
this in the penultimate paragraph of the Introduction is
satisfied.\qfd

\begin{lema}
The set of neighbors of each copy
$\zeta_n^{i,j}(\vec{ST}_n)=\vec{ST}_n^{i,j}$ of $\vec{ST}_n$ in
$\vec{ST}_{n+1}$ is an {\rm E}$_\pm$-set that induces an
E$_\pm$-subdigraph $ST_n^i$ in $\vec{ST}_{n+1}$ $($independen\-tly
 in $j)$ isomorphic to the star graph $ST_n$ considered oriented from
the vertices with $0$th entry equal to $i$ into the vertices with
$1$th entry equal to $i$.
\end{lema}

\proof For each $n>1$\,, $D_n=\vec{ST}_{n+1}$ is devised so as to be
member of a neighborly chain, via the inclusive map
$\zeta_n^{n,n}=\kappa_n^{n,n}$. By permuting coordinates, it is seen
that a similar behavior occurs for any other map $\zeta_n^{i,j}$. In
fact, each vertex $v$ in a copy $\vec{ST}_n^{i,j}$ of $\vec{ST}_n$
in $\vec{ST}_{n+1}$ is the only intersecting element of such
$\vec{ST}_n^{i,j}$ with a specific directed triangle
$\vec{\Delta}_v$. (For example, vertex $v=13042$ of $\vec{ST}_5$ in
Figure 3 is the only intersecting vertex of $\vec{ST}_4^{0,2}$ with
the directed triangle $\vec{\Delta}_v$ whose other two vertices are
$u=01342$ and $w=30142$). Then, the subdigraph of $\vec{ST}_{n+1}$
induced by all subdigraphs $\vec{\Delta}_v-\{v\}$\,, where $v\in
V(\vec{ST}_n^{i,j})$\,, (each such $\vec{\Delta}_v-\{v\}$ containing
just one arc), is the claimed E$_\pm$-subdigraph $ST_n^i$\,,
independently of $j$\,, and is clearly isomorphic to $ST_n$\,,
oriented as stated. \qfd

\begin{lema} Consider the E$_\pm$-subdigraph $ST_n^i$ of {\rm Lemma 4}.
There are\, $n+1$\, copies $ST_n^i$ of the star graph $ST_n$ in
$\vec{ST}_{n+1}$\,, where $0\le i\le n$\,, with their arcs oriented
from their vertices with $0$th entry equal to $i$ onto their
vertices with $1$th entry equal to $i$. The set of neighbors of each
of these copies $ST_n^i$ of $ST_n$ is the disjoint union of $n-1$
copies of $\vec{ST}_n$ counting once the image of each map
$\zeta_n^{i,j}$\,, where $2\le j\le n$.
\end{lema}

\proof Clearly, the index $i$ in $ST_n^i$ varies in
$\{0,1,\ldots,n\}$\,, yielding the $n+1$ copies of $ST_n^i$ in the
statement. Each such copy is induced by the set of $(+)$neighbors
and $(-)$neighbors of the vertices of any fixed $\vec{ST}_n^{i,j}$.
Since the index $j$ here may be any value in $\{2,\ldots,n\}$\,, the
second sentence of the statement holds, too.\qfd\bigskip

\noindent As a case to exemplify the condition of neighborly
E$_\pm$-chain in the proof of Theorem 6 below, the assignment $\rho$
from $V(\vec{ST}_n^{n,n})=V(\kappa_n^{n,n}(\vec{ST}_n))$\,, for
$n=3$\,, onto the disjoint union $K$ of $|V(\vec{ST}_n^{n,n}))|=3$
digraphs $\vec{P_2}$ indicated in the penultimate paragraph of the
Introduction, is given by:
$$^{0123\mapsto(3021,1320),\;\; 1203\mapsto(3102,2301),\;\;
2013->(3210,0312),}$$ (represented on the right quarter of Figure
2). Here, each assignment $v\mapsto\rho(v)$ yields a corresponding
directed triangle $\vec{\Delta}_v$\,, induced by $v$ and
$\rho(v)$\,, namely

\noindent $(0123,3021,1320)$, $(1203,3102,2301)$ and $(2013,3210,0312)$,
respectively, where the 0th, 1th and 3th coordinates are modified in
the 3 triangles as specified and the remaining coordinate remains
fixed (in the respective values 2, 0 and 1). Replacing $n=3$ by
$n=4$\,, the assignment $\rho$ is given now by:
$$\begin{array}{c}^{
32104\mapsto(43102,24103);\;\; 20134\mapsto(42130,04132);\;\;
21304\mapsto(42301,14302);}_{ 01234\mapsto(40231,14230);\;\;
13204\mapsto(41203,34201);\;\;12034\mapsto(41032,24031);}\vspace*{0.5mm}\\
^{02314\mapsto(40312,24310);\;\; 23014\mapsto(42013,34012);\;\;
31024\mapsto(43021,14023);}_{ 30214\mapsto(43210,04213);\;\;
10324\mapsto(41320,04321);\;\; 03124\mapsto(40123,34120).}
\end{array}$$
where the sources in the 12 cases are shown as those common to two
light-gray equilateral triangles in Figure 1 from left to right and
from top to bottom.

\begin{theorem}
The star digraphs $\vec{ST}_n$\,, where $n\ge 1$\,, are strong and
constitute a dense segmental neighborly {\rm E}$_\pm$-chain.
\end{theorem}

\proof For $n>2$\,, Lemma 4 insures that $\vec{ST}_n$ contains a
copy of $ST_{n-1}$\,, induced by an E$_\pm$-set. The undirected
version of this copy contains a Hamilton cycle $H$ \cite{CKRR,KL}
and allows to show that $\vec{ST}_n$ is strong: Given vertices $u,v$
in $\vec{ST}_n$\,, there is a path $P=uw_0\ldots w_kv$ in
$\vec{ST}_n$\,, where $w_0\ldots w_k$ is a section of $H$; $P$ is
transformed into a directed path by replacing each backward arc
$\overleftarrow{a}$ of it by the directed 2-path forming an oriented
triangle of $\vec{ST}_n$ with it. Again by Lemma 4, the star
digraphs form an E$_\pm$-chain $\mathcal D$ as expressed in (5).
Since they satisfy equality (1), $\mathcal D$ is dense. The
inclusive maps $\kappa_{n-1}=\kappa_n^{n,n}=\zeta_n^{n,n}$ from
$D_{n-1}=\vec{ST}_n$ into $D_n=\vec{ST}_{n+1}$\,, ($n>1$), show that
$\mathcal D$ is neighborly, because
$S_n=N^+(V(\kappa_{n-1}(D_{n-1})))\cup
N^-(V(\kappa_{n-1}(D_{n-1})))$ is a disjoint union of two stable
vertex subsets of $D_n$ as indicated and there is a bijective
correspondence $\rho:V(\kappa_n(D_{n-1}))\rightarrow K$ such that
$v$ and $\rho(v)$ induce a directed triangle $\vec{\Delta}_v$\,, for
each vertex $v$ of $\kappa_{n-1}(D_{n-1})$, where $K$ is a disjoint
union of $|V(\kappa_{n-1}(D_{n-1}))|$ digraphs $\vec{P}_2$ in $D$
consisting each of a single arc from $N^+(V(\kappa_{n-1}(D_{n-1})))$
to $N^-(V(\kappa_{n-1}(D_{n-1})))$. In order to establish that the
star digraphs $\vec{ST}_n$ form a segmental E$_\pm$-chain, the
examples of partitions in the last paragraph of Section 4 can now be
directly generalized. For example, the E$_\pm$-set $ST_n^0$ of
$\vec{ST}_{n+1}$ admits $n-1$ different partitions into $n$ copies
of $ST_{n-1}$, namely
$\{\zeta_n^{i,j}(ST_{n-1}^0)\,|\,i=1,\ldots,n\}$, for
$j=n,\ldots,3,2$. In these partitions, each copy of $ST_{n-1}$, like
the 6-cycle dags for $n=3$ in the table of the mentioned paragraph,
is the induced subdigraph of a E$_\pm$-set $S_n$ in a corresponding
subdigraph $D_n=\vec{ST}_{n+1}^{i,j}$, which is the requirement for
every $n\ge 1$ cited in that paragraph in order to insure that
$\mathcal D$ is inclusive. Since partitions as in the mentioned
table are obtained for every $n\ge 1$, we conclude that $\mathcal D$
is segmental.
 \qfd

\section{Some comments and open problems}

\subsection{Hamiltonicity and traceability}

\begin{ques} Are all the star digraphs traceable? Hamiltonian?
\end{ques}

\noindent The star digraph $\vec{ST}_4$ is not hamiltonian. We think
that this is the case for every star digraph $\vec{ST}_n$\,, $n>3$.
For example, there are just two types of oriented Hamilton paths in
$\vec{ST}_4$\,, obtained as follows. Let us start a path $P$ at a
fixed vertex $v$ of $\vec{ST}_4$\,, indicating by $b$ whenever a
2-arc is added to $P$\,, and by $a$ whenever just a 1-arc is added
to $P$\,, (steps represented respectively by two subsequent arcs and
by just one arc bordering an directed triangle $\vec{\Delta}$ in
$\vec{ST}_4$); then we get the claimed two types: $P=aababbb$\,,
obtained by setting a starting $a$ and having continuation
preference for $a$ over $b$ unless backtracking is necessary in
trying to produce a Hamilton path, and the reversal
$P^{-1}=bbbabaa$.

\subsection{Pancake digraphs}

\noindent For $n>4$\,, the {\it pancake digraph} $\vec{PC}_n$ is
defined as the oriented Cayley graph of $Sym_n$ with respect to the
set of compositions $(0\;1)\circ f$\,, where $f$ runs over the set
of involutions $\{\Pi_{j=1}^{\lfloor
i/2\rfloor}(j\;(i-j))=(1\;i)(2\;(i-1))\cdots(\lfloor i/2\rfloor\;
\lceil i/2 \rceil);\; i\in I_n\setminus\{0,1\}\}$. Such a
$\vec{PC}_n$ is connected but its definition could not hold for
$n\le 4$ if we are to keep connectedness, as such a digraph would
have two components, both isomorphic to $\vec{ST}_n$. In particular,
$\vec{PC}_5$ is obtained from two disjoint copies of $\vec{ST}_5$
(one with vertex set $Alt_5$\,, the other with vertex set
$(Sym_5\!\setminus\! Alt_5)$), and replacing the pairs of arcs
corresponding to the right multiplication by the generator
$(01)(14)(23)=(041)(23)$ of $\vec{PC}_5$ (as a Cayley graph) by
corresponding crossed arcs between the two said copies of
$\vec{ST}_5$. In fact, we could maintain the toroidal cutout of
Figures 2 and 3 while replacing each vertex $a_0a_1a_2a_3a_4$ by
$a_0a_1a_3a_2a_4$ (permuting $a_2$ and $a_3$) in order to obtain a
copy of $\vec{ST}_5$ with vertices replaced from $Alt_5$ to
$(Sym_5\!\setminus\! Alt_5)$. Let us call this second copy of
$\vec{ST}_5$ by $\overleftarrow{ST}_5$. Then, $\vec{PC}_5$ is
obtained by modifying the disjoint union
$\vec{ST}_5\cup\overleftarrow{ST}_5$ by replacing each pair of arcs
$\{(a_0a_1a_2a_3a_4,a_4a_0a_2a_3a_1),(a_0a_1a_3a_2a_4,a_4a_0a_3a_2a_1)\}$
by the pair of crossed arcs
$\{(a_0a_1a_2a_3a_4,a_4a_0a_3a_2a_1),(a_0a_1a_3a_2a_4,a_4a_0a_2a_3a_1)\}.$
In $\vec{ST}_5$\,, the arcs of the form
$(a_0a_1a_2a_3a_4,a_4a_0a_3a_4a_1)$ induce the disjoint union of 20
directed triangles which are the intersections of the pairs of
copies of the form $\{\vec{ST}_4^{i,2},\vec{ST}_4^{j,3};i\ne j;\}$.
In the way from $\vec{ST}_5$ to $\vec{PC}_5$\,, these 20 triangles
give place to 20 corresponding oriented 6-cycles, each formed by 3
pairs of crossed pairs as above. In addition, the 40 remaining
directed triangles of $\vec{ST}_5$ give place to a total of 80
directed triangles in $\vec{PC}_5$. We recall from \cite{io} that
the pancake graphs $PC_n$ form a dense segmental neighborly E-chain.

\begin{ques} For $n>4$\,, do the pancake digraphs $\vec{PC}_n$ form a dense
segmental neighborly {\rm E}$_\pm$-chain? Are they strong?
Traceable? Hamiltonian?
\end{ques}

\subsection{Binary-star digraphs}

\noindent A different variant of the star digraphs $\vec{ST}_n$\,,
on $n!$ vertices (like the pancake graph $\vec{PC}_n$) is the {\it
binary-star digraphs} $B\vec{ST}_n$\,, defined as the bipartite
graph whose vertex parts are the cosets of $Alt_n$ in $Sym_n$\,,
with an arc $(\sigma,\sigma\circ(1\;i))$ for each $\sigma\in
Alt_n$\,, and an arc $(\sigma,\sigma\circ(0\;i))$\,, for each
$\sigma\in (Sym_n\!\setminus\! Alt_n)$\,, where $i\in
I_n\setminus\sigma \{0,1\}$. The reader is invited to check that
 $B\vec{ST}_n$ is isomorphic to
the canonical $2$-covering bipartite digraph of $\vec{ST}_n$.

\begin{ques}
Do the binary-star digraphs form a dense segmental neighborly
{\rm E}$_\pm$-chain?
Are they strong? Traceable?
Hamiltonian? Strongly Hamiltonian traceable as in
{\rm\cite{CKRR,hch}}, in a directed sense? Hamiltonian connected, as
conjectured in {\rm\cite{CKRR,HS}} for the star graphs?\end{ques}

\begin{ques} Do there exist infinite families of {\rm E}$_\pm$-chains of
Cayley digraphs on symmetric groups that include both the
binary-star and pancake digraphs, in a fashion similar to Section
$2$ of {\rm\cite{io}}?
\end{ques}

\end{document}